\documentclass[12pt]{article}
\usepackage{amsmath}[1996/11/01]
\usepackage{amssymb,amsthm,amsfonts,latexsym,psfig,epsfig,graphics}
\usepackage{amscd}

\advance \textheight by -1 \topmargin
\advance \textheight by -1 \headheight
\advance \textheight by -1 \headsep
\advance \textheight by -2 \footskip
\vsize = \textheight
\hsize = \textwidth
\newtheorem{theorem}{Theorem}[section]

\newtheorem{kor}[theorem]{Corollary}

\newtheorem{rem}[theorem]{Remark}
\newtheorem{properties}[theorem]{Properties}
\newtheorem{lemma}[theorem]{Lemma}

\newtheorem{definition}[theorem]{Definition}

\renewcommand{\setminus}{-}

\newcommand{\be}{\begin{enumerate}}
\newcommand{\ee}{\end{enumerate}}
\newcommand{\bi}{\begin{itemize}}
\newcommand{\ei}{\end{itemize}}

\newcommand{\ba}{\begin{array}}
\newcommand{\ea}{\end{array}}

\newcommand{\Z}{{\mathbb{Z}}}
\newcommand{\Q}{{\mathbb{Q}}}
\newcommand{\F}{{\mathbb{F}}}
\newcommand{\N}{{\mathbb{N}}}

\newcommand{\tr}{\mbox{trace}}

\newcommand{\eb}{\phantom{zzz}\hfill{$\square $}\smallskip}

\renewcommand{\em}{\sf}

\DeclareMathOperator{\End}{End}
\DeclareMathOperator{\Ext}{Ext}
\DeclareMathOperator{\diag}{diag}
\DeclareMathOperator{\Id}{Id}
\DeclareMathOperator{\Tr}{Tr}
\DeclareMathOperator{\Gal}{Gal}

\begin{document}

\Huge
\begin{center}
{\bf On blocks with cyclic defect group and their head orders.}
\end{center}
\normalsize
\begin{center}
Gabriele Nebe \footnote{Abteilung Reine Mathematik, Universit\"at Ulm, 89069 Ulm, Germany, nebe@mathematik.uni-ulm.de}
\footnote{Radcliffe fellow}
\end{center}

\small

{\sc Abstract}:
It is shown that 
\cite[Theorem 8.5]{Ple 83} 
describes blocks of cyclic defect group up to Morita equivalence.
In particular such a block is determined %(up to Morita equivalence)
by its planar embedded Brauer tree. % (with permutations and multiplicities).
Applying the radical idealizer process the head order of such blocks
is calculated explicitly.

\normalsize 

\section{Introduction.}
Blocks with cyclic defect group are very well understood.
Despite of their very special structure
these blocks are extensively used to study examples for the validity of
various conjectures, since they are essentially described by combinatorial
means.
A detailed introduction to the theory of blocks with cyclic defect groups,
that also deals with rationality questions of the involved characters, 
is given in \cite[Chapter 7]{Feit}.
Using the known character theoretic 
 information and some new methods, essentially based on linear algebra,
Plesken \cite[Chapter 8]{Ple 83} gives a rather explicit description of
blocks 
 with cyclic defect group
$ B = \Z_pG\epsilon$
 of $p$-adic group rings.
Theorem \ref{unique} shows
 that \cite[Theorem 8.5,8.10]{Ple 83} determines the ring theoretic structure
of such blocks up to isomorphism.
In particular the planar embedded Brauer tree
together with the character fields and the Galois action on the characters 
determine the Morita equivalence class of $B$.

The second part of this paper deals with the {\em radical idealiser chain} of
$ B$. This is a finite chain associated to an order $\Lambda $
$$\Lambda  =: \Lambda _0 \subset \Lambda _1 \subset \ldots \subset \Lambda _N = \Lambda _{N+1}$$
where $\Lambda _{i+1} = \Id (J(\Lambda _{i})) $ $(i=0,\ldots , N)$
is the 2-sided  idealizer of the Jacobson radical of $\Lambda _i$
that necessarily ends in a hereditary order $\Lambda _N$ called the 
{\em head order} of $\Lambda $.
Section \ref{head} investigates the radical idealiser chain of blocks 
with cyclic defect group and 
calculates the head order of such a block.

\section{Blocks with cyclic defect group.}

The aim of this section is to show, how far \cite[Theorem 8.5]{Ple 83} 
determines blocks with cyclic defect group  of group rings over discrete
valuation rings. Moreover we provide the tiny extra bit of information to
get a complete description up to isomorphism.

Throughout the paper let $R$ be a (not necessarily commutative)
discrete valuation ring with prime element $\pi $
and residue class field $k = R/\pi R$ and let $K$ be the 
skew-field of fractions of $R$.

A convenient language to describe certain $R$-orders are exponent matrices.

\begin{definition} (see \cite[Definition 39.2]{Rei}, \cite{Ple 83})
Let $d=(d_1,\ldots, d_n ) \in \N ^n $, $D:=\sum _{i=1}^n d_i $ and
$M\in \Z^{n\times n }$.
Then 
$$\Lambda (R,d,M) := \{ X = (x_{ij}) \in K^{D\times D } \mid x_{ij} \in \pi ^{m_{ij}} R^{d_i\times d_j}  \} .$$
\end{definition}

\underline{Example:} (hereditary orders, see \cite[Section 39]{Rei})
Let $\Lambda $ be a hereditary order and $P$ be a projective $\Lambda $-lattice.
Then $\overline{P} := P/\pi P$ is uniserial. Let $d_1,\ldots , d_n$ be the dimensions of the
simple $\Lambda $-modules in the order in which they occur in the radical series
of $\overline{P}$.
Then with respect to a suitable $R$-basis of $P$ (adapted to this lattice
chain), 
$$\Lambda = \Lambda (R,d , H_n) , \mbox{ where } 
H_n = \left( \begin{array}{ccccc}
0 & 1 & \ldots & \ldots & 1 \\ 
0 & 0 & 1 & \ldots & 1  \\
\vdots & \ddots & \ddots & \ddots & \vdots \\
0 & \ldots & \ldots & 0 & 1 \\
0 & \ldots & \ldots & 0 & 0 \end{array} \right) $$

The description of blocks with cyclic defect  given in \cite{Ple 83} will be
repeated briefly.
All the following results can be found in Chapter 8 of this lecture notes, so we 
omit the detailed citations.
Let $K$ be an {\bf unramified} extension of $\Q _p$, $R$ its ring of integers and
${\cal B} $ be a block of $RG$ with cyclic defect group of order $p^a$.
Assume that $k:= R/pR $ is a splitting field of $k {\cal B} $.
By \cite[Chapter 7]{Feit}, the minimal choice of such a field $K$ 
is the character field of any of the non-exceptional characters in ${\cal B}$.

Let $\epsilon _1,\ldots , \epsilon _h$ be the central primitive idempotents
in ${\cal A}:= K {\cal B} $.
Then $h=a+e$ where $e$ is the number of simple $k{\cal B} $-modules and 
after a suitable permutation of the indices 
the centers $Z({\cal A}\epsilon _s) \cong K$ for $s=a+1,\ldots , a+e$ and 
$Z({\cal A}\epsilon _s) =: Z_s$ is a totally ramified extension of $K$ of degree
$\frac{p^s-p^{s-1}}{e}$ for $s=1,\ldots , a$.
The center of ${\cal B} \epsilon_s $ is the maximal order $R_s$ in $Z_s$.

The vertices in the Brauer tree are the exceptional vertex $\{ 1,\ldots , a \}$ 
and $a+1,\ldots , a+e $ corresponding to the other simple ${\cal A}$-modules.
Let $T_{odd}$ resp. $T_{even}$ 
 denote the set of vertices having an odd (resp. even) distance from the exceptional
vertex.  

For $s\in \{1,\ldots , h \}$ let $r_s \subset \{ 1,\ldots , e\}$ be the set of indices of the 
simple constituents of any ${\cal B} \epsilon _s $ lattice.
Then for the exceptional vertex $r_1= \ldots = r_a $ and 
the sets $r_s$ are the orbits of certain permutations $\delta  $ (if $s\in T_{even}$) 
resp. $\rho $ (if $s\in T_{odd}$).
Let $d_1,\ldots , d_e $ be the $k$-dimensions of the simple $k{\cal B} $-modules
and $f_1,\ldots , f_e $ be orthogonal idempotents of ${\cal B} $ that lift 
the corresponding central primitive idempotents of ${\cal B} / J({\cal B} )$.

\begin{theorem}\label{cycple1} (\cite[Theorem 8.3]{Ple 83})
With the notation above let $i\in r_s$. Then
$$(i) \ \ \ \ {\cal B} \epsilon _s \cong \Lambda (R_s, (d_i, d_{\delta (i) } ,\ldots , d_{\delta ^{|r_s|-1}(i)} ), H_{|r_s|} ) \
\mbox{ for } s=1,\ldots , a $$ 
and 
$$(ii) \ \ \ \ {\cal B} \epsilon _s \cong \Lambda (R, (d_i, d_{\sigma (i) } ,\ldots , d_{\sigma ^{|r_s|-1}(i)}) , a H_{|r_s|} ) \ 
\mbox{ for } s=a+1,\ldots , a+e $$ 
where $\sigma = \delta $ if $s\in T_{even}$ and $\sigma = \rho $ if $s\in T_{odd}$.
\end{theorem}

It remains to describe how ${\cal B} $ sits inside the direct sum of 
the ${\cal B} \epsilon _s$, that is to describe the amalgamations between the
${\cal B} \epsilon _s$.

\begin{theorem} \label{cycple} (\cite[Theorem 8.5]{Ple 83})
\begin{itemize}
\item[(i)]
For the exceptional vertex 
$\Gamma _{a} := (\epsilon _1 + \ldots + \epsilon _a) {\cal B} $ one gets an 
inductive description:
For $s=2,\ldots , a$ let 
$$X_s := {\cal B} \epsilon _s / J({\cal B} \epsilon _s)^{x_s} =
{\cal B} \epsilon _s / \pi _s^{y_s} {\cal B} \epsilon _s $$
where $x_s = |r_1| y_s $ and $y_s = \frac{p^{s-1}-1}{e} $ and
let $\nu _s : {\cal B} \epsilon _s \to X_s $ be the natural epimorphism.
Define $R$-orders $\Gamma _s$ $(s=1,\ldots , a )$ inductively 
by $$\Gamma _1 := \epsilon _1 {\cal B}  \mbox{ and  }
\Gamma _{s} := \{ (x,y) \in \Gamma _{s-1} \oplus {\cal B}  \epsilon _s \mid 
\varphi _{s-1} (x) = \nu _s (y) \} $$
where $\varphi _{s-1} $ is an epimorphism from $\Gamma _{s-1} $ onto $X_s$.
\item[(ii)]
Let $\Gamma _0 := (\epsilon _{a+1} + \ldots + \epsilon _{a+e}) {\cal B} $.
Then for any $i,j \in \{ 1,\ldots , e \} \setminus r_1 $ with $i\neq j$ one 
gets 
$$f_i \Gamma _0 f_j = \oplus _{s=a+1}^{a+e} f_i {\cal B} \epsilon _s f_ j $$
and 
$$f_i \Gamma _0 f_i \cong \{ (x,y) \mid x,y \in R^{d_i \times d_i } , 
x \equiv y \pmod{p^a} \} \subset f_i \epsilon _s {\cal B} f_i \oplus f_i \epsilon _t {\cal B} f_i \cong R^{d_i\times d_i} \oplus R^{d_i \times d_i }$$
if $i\in r_s \cap r_t $.
\item[(iii)]
Finally there are epimorphisms 
$\nu $ and $\mu$ of $\Gamma _0$ and $\Gamma _a$ onto 
$\oplus _{i\in r_1} (R/p^a R)^{d_i \times d_i }$ such that 
$${\cal B} = \{ (x,y) \in \Gamma _0 \oplus \Gamma _a \mid \nu (x) = \mu (y) \} .$$
\end{itemize}
\end{theorem}

Note that 
$X_s \cong  \Gamma _{s-1} / J(\Gamma _{s-1})^{x_s} $.

The possible ambiguity in this description is the choice of the epimorphisms in
(i) and in (iii).
It is clear that one can always fix one of the two epimorphisms.
The choices for the other one correspond to the automorphisms of the image.
So the question is, whether these automorphisms lift to automorphisms of 
${\cal B} $. This is clear for the maps in (iii).
For (i) this is unfortunately not always the case. 

To simplify notation, it is convenient to pass to the Morita-equivalent
basic order. 
Let $S$ be any discrete valuation ring with prime element $\pi _S $ and
$\Lambda := \Lambda (S,(1,\ldots , 1), H_n) $ the basic hereditary 
$S$-order of degree $n$.
Let $X:= \Lambda / \pi  _S\Lambda $.
Then $\Lambda $ is generated by the idempotents 
$e_i = \diag (0,\ldots, 0,1,0,\ldots 0) $ (the $1$ is on the 
$i$-th place), the elements 
$$e_{i+1,i} = \left( \begin{array}{cccc}
0 & \ldots & \ldots & 0 \\
\vdots & \ddots & \ddots & \vdots \\
\vdots & 1 & \ddots & \vdots \\
0 & \ldots & \ldots & 0 \end{array} \right) \in e_{i+1} \Lambda e_i  
\mbox{ and }
g_{1,n}  = \
 \left( \begin{array}{cccc}
0 & \ldots & 0 & \pi  _S\\
0 & \ldots & \ldots & 0 \\
\vdots & \ddots & \ddots & \vdots \\
0 & \ldots & \ldots & 0 \end{array} \right) \in e_{1} \Lambda e_n  $$
where $i=1,\ldots , n-1$.
These generators map onto generators $\overline{e_i}$, $\overline{e_{i+1,1}}$
and $\overline{g_{1,n}}$ of $X$ corresponding to the
Ext-quiver of $X$ which is a directed $n$-gon.
They satisfy the relation that 
$$\overline{g_{1,n}} \overline{e_{n,n-1}} \cdots \overline{e_{2,1}} = 0 $$
and similarly for any cyclic permutation of this product.

\begin{lemma}\label{basechange}
Let $\Lambda $ and $X = \Lambda/\pi  _S\Lambda $ be as above and let
$\varphi $ be an automorphism of $X$ that fixes all the idempotents
$\overline{e_i}$.

Then there are $0 \neq \overline{\lambda _i} \in  S/\pi  _SS =: k _S$ with
$\varphi ( \overline{e_{i+1,i}} ) = \overline{\lambda _i} \overline{e_{i+1,i}}  $ for 
$i=1,\ldots , n-1$ and
$\varphi ( \overline{g_{1,n}} ) = \overline{\lambda _n} \overline{g_{1,n}}  $.

There is an automorphism $\phi $ of $\Lambda $ that lifts $\varphi $
if and only if the product $\overline{\lambda_1}\cdots \overline{\lambda _n} = 1$.

In particular, there is always an  automorphism $\phi $ of $\Lambda $ with
$$\overline{\phi (e_{i+1,i}) }  = \varphi (\overline{e_{i+1,i}}) , \mbox{ and }
\phi (e_j ) = e_j \mbox{ for all } i=1,\ldots, n-1, j=1,\ldots , n .$$
\end{lemma}

\proof
The automorphism $\varphi $ maps the generator 
$\overline{e_{i+1,i}} \in \overline{e_{i+1}} X \overline{e_i} = k _S\overline{e_{i+1,i}}$ to 
some other generator of this module
($i=1,\ldots , n-1$) and
similar for $\overline{g_{1,n}}$.
Hence there are such units $\overline{\lambda _i} \in k _S^*$ as described in the lemma.
Moreover any such tuple $( \overline{\lambda _1},\ldots , \overline{\lambda _n})
\in (k_S^*)^n$
determines a unique automorphism of $X$ fixing all the idempotents $\overline{e_i}$.

Choose units $\lambda _i \in S^*$ that map to 
$\overline{\lambda _i} $ in $k _S$.
Then the matrix $$D:=\diag(1,\lambda _1 ,\lambda _1\lambda _2 ,\ldots , \lambda _1 \cdots \lambda _{n-1} ) \in \Lambda ^*$$
fixes all the $e_i$ and conjugates $e_{i+1,i} $ to 
$\lambda _i e_{i+1,i}$ for all $i=1,\ldots , n-1 $ and 
$g_{1,n} $ to $(\lambda _1 \cdots \lambda _{n-1})^{-1} g_{1,n}$.
Hence if the product of the $\overline{\lambda _i}$ is 1, then conjugation
by $D$ is the desired automorphism $\phi $.

On the other hand it is easy to see that all automorphisms of $\Lambda $ that
fix the idempotents $e_i$ are given by conjugation with a diagonal matrix
$D=\diag(d_1,\ldots , d_n)$ mapping the matrix units $e_{ij}$ to $\frac{d_i}{d_j} e_{ij}$.
\eb

Therefore there is a tiny bit missing in Theorem \ref{cycple} (i) to describe the
exceptional vertex $\Gamma _a$ up to isomorphism.
Since blocks of group rings are symmetric orders, however, the missing information
can easily be obtained from the trace bilinear form. 

Theorem \ref{cycple} gives the center $Z:=Z(\Gamma _a)$ up to isomorphism.
Instead of continuing with Plesken's description, it seems to be easier to give
generators for $\Gamma _a$ over the center $Z$ using the $\Ext $-quiver of
$\Gamma _a$.
To this aim, we pass to the Morita equivalent basic order and assume that
all simple $\Gamma _a$-modules are of dimension 1 over $k$.
All projective $\Gamma _a$-lattices  are uniserial when reduced modulo $p$,
 where the
sequence of composition factors is given by the permutation $\delta $.
Therefore the $\Ext $-quiver of $\Gamma _a$ is again a directed $n$-gon, where
$n= |r_1|$ is the number of simple $\Gamma _a$-modules.
If $n=1$, then $\Gamma _a = Z(\Gamma _a)$ is already described completely by
Theorem \ref{cycple}.
So we will assume that $n>1$.
Let $e_1,\ldots,e_n \in \Gamma _a$ be orthogonal lifts of the central primitive idempotents of 
$\Gamma _a/ J(\Gamma _a)$ ordered in such a way that 
$\delta _{|r_1} = (1,\ldots , n)$.
Denote the corresponding matrix units 
in $e_i K\Gamma _a  \epsilon _s e_j $ by $e_{ij}^s$ ($i,j\in \{1,\ldots , n\}, s\in \{1,\ldots , a \}$).
Then according to Theorem \ref{cycple} and Lemma \ref{basechange}
after a choice of a suitable basis generators of $\Gamma _a$ over its center $Z$ 
can be chosen as 
$$e_1,\ldots , e_n, (e_{2,1}^1,\ldots , e_{2,1}^a) =: e_{2,1} \in
e_2 \Gamma _a e_1 , \ldots , 
(e_{n,n-1}^1 ,\ldots , e_{n,n-1}^a) =: e_{n,n-1} \in e_n \Gamma _a e_{n-1} $$ and 
$$ 
(x_1\pi_1 e_{1,n}^1 ,\ldots , x_a\pi_a e_{1,n}^a) =: g_{1,n} \in 
e_1 \Gamma _a e_n
\mbox{ for certain units } x_i \in Z_i = Z(\Gamma _a \epsilon _i). $$

\begin{theorem}\label{unique}
Let $Z:=Z(\Gamma _a)$ and let $Z^{\# }$ be the dual of $Z$ with respect to the
sum of usual the trace bilinear forms.
Then there are units $x_i\in Z_i = Z(\Gamma _a \epsilon _i)$ ($i=1,\ldots , a$)
such that 
$$p^a Z ^{\# } = (x_1 \pi _1,\ldots , x_a \pi_a ) Z .$$
With the choice of these $x_i$, the order $\Gamma _a$ is generated by
$$Z,e_1,\ldots , e_n, e_{i+1,i} \ (i=1,\ldots , n-1), \mbox{ and } g_{1,n} $$
as defined above.
\end{theorem}

\proof
We may assume that $n>1$.
${\cal B}$ is a symmetric order with respect to the 
associative symmetric bilinear form
$$(x,y) \mapsto 
 \frac{1}{|G|} \tr _{reg} (x y ) = \tr_{red}(x y z) = : \Tr _z(x,y) $$
where $\tr _{reg}$ and $\tr _{red} $ denote the regular 
respectively reduced trace of $K{\cal B} $
and $z = \sum _{s=1}^{a+e} \frac{\chi _s(1)}{|G|} \epsilon _s $,
where $\epsilon _1,\ldots , \epsilon _{a+e}$ are the central primitive idempotents
of $K{\cal B} $ and $\chi _1,\ldots , \chi_{a+e}$ corresponding absolutely
irreducible (complex) characters of $G$.

Let $f_1,\ldots , f_n$ denote orthogonal idempotents in ${\cal B}$ that map onto the
central primitive idempotents of ${\cal B}/J({\cal B})$ such that $$e_i = f_i (\epsilon _1+\ldots + \epsilon _a ) \ (i=1,\ldots , n ).$$
Since $n>1$ 
$$\langle g_{1,n} \rangle _Z = 
e_1 \Gamma _a e_n = f_1 {\cal B} f_n  = (f_n {\cal B} f_1 )^{\#} = (e_n \Gamma _a e_1 )^{\#} =
\langle e_{n,n-1} \cdots e_{2,1} \rangle _Z ^{\#} $$
can be calculated via the symmetrizing form above.
Since the character degrees of the absolutely irreducible characters
belonging to the exceptional vertex are all equal
%$$\chi _1(1) = \ldots = \chi _a (1)  ,$$
the dual with respect to $\Tr _z$ is as stated in the theorem, 
yielding the remaining generator $g_{1,n}$ for $\Gamma _a$.
\eb

Note that the $x_i$ do not depend on the degrees of the 
irreducible complex characters in ${\cal B}$, since all exceptional
absolutely irreducible characters have the same degree.
Therefore one gets 

\begin{kor}{\label{unram}}
Let ${\cal B}_i$ ($i=1,2$) be two blocks with
cyclic defect group $\cong C_{p^a}$ and
assume that $R$ is an unramified extension of $\Z_p$ 
that is large enough so that $k$ is a splitting field for 
$k{\cal B}_i$.
Then ${\cal B}_1$ and ${\cal B}_2$ are Morita equivalent 
if and only if their Brauer trees (including the permutations $\delta $ and $\rho $)
and the 
character fields $Z_1 ,\ldots , Z_a$ coincide.
\end{kor}

Also, symmetric orders remain symmetric orders, when one extends the 
ground ring. 
Therefore 
the  explicit description in \cite[Theorem 8.5]{Ple 83} shows that
the Brauer tree determines a block of cyclic defect up 
to Morita equivalence (over an algebraically closed field).
This is also shown in \cite[Theorem 2.7(ii)]{Lin96} with completely
different methods.

\begin{kor}{\label{Morita}}
Let ${\cal B}_i$ ($i=1,2$) be two blocks with isomorphic
cyclic defect group and
assume that $R$ is large enough so that $k$ and $K$ are splitting fields for 
$k{\cal B}_i$ and $K{\cal B}_i$. (Here we drop the assumption that 
$K$ is unramified over $\Q_p$.)
Then ${\cal B}_1$ and ${\cal B}_2$ are Morita equivalent 
if and only if their planar embedded Brauer trees coincide.
\end{kor}

\subsection{Galois descent.}
We now perform the Galois descent to obtain a description over $\Z_p$
(see \cite[Chapter 8]{Ple 83}).
So let $B$ be a block of $\Z_pG$ such that ${\cal B}$ is a summand of 
$R\otimes B$. 
We assume that $K$ is chosen to be minimal, i.e. 
$K = \Q _p[\chi_{a+1}] = \ldots =\Q _p[\chi _{a+e}]$ 
is the character field of any non-exceptional
absolutely irreducible Frobenius character that belongs to ${\cal B}$.
The maximal unramified subfield $\tilde{K}$ of the character field 
$\tilde{Z}_s := \Q _p [\chi _s]$ $(s=1,\ldots , a )$ of any 
exceptional absolutely irreducible Frobenius character in ${\cal B}$ 
does not depend on the character and is a subfield of $K$.
Let $m:=[K:\tilde{K}]$ denote its index.

If $\tilde{R}$ denotes the ring of integers in $\tilde{K}$, then
$\tilde{R}$ embeds into the center of $B$ such that  $B$ 
can be viewed as an $\tilde{R}$-order and
$R\otimes _{\tilde{R}} B \cong {\cal B}.$

The Galois group $\Gal (K/\tilde{K})  = \Gal (k/\tilde{k}) \cong C_m$
(where $\tilde{k} := \tilde{R}/p\tilde{R}$) acts on the simple 
${\cal B}$-modules and the corresponding idempotents $f_1,\ldots , f_e$
with orbits of length $m$.
Therefore orthogonal lifts of the central primitive idempotents of
$B/J(B)$ can
be chosen as $\tilde{f}_1,\ldots , \tilde{f} _{\tilde{e}} \in B$ 
where $\tilde{e}:=\frac{e}{m}$
is the number of simple 
$\F_pB $-modules, each of which has character field $k = R/pR$.

The central primitive idempotents in $A:=\Q_p \otimes B$ are 
$\tilde{\epsilon }_1,\ldots , \tilde{\epsilon }_a ,
\tilde{\epsilon }_{a+1},\ldots , \tilde{\epsilon }_{a+\tilde{e}}  $
indexed in such a way that
$\tilde{Z}_s$ is a totally 
ramified extension of $\tilde{K}$ of degree
$\frac{p^s-p^{s-1}}{e}$ for $s=1,\ldots , a$.

For an appropriate 
ordering of the index set $\{1,\ldots , e\}$ 
the $k$-dimensions of the simple $\F_pB$-modules are 
$d_1,\ldots , d_{\tilde{e}}$
and 
the set of  indices of the simple $\tilde{\epsilon }_s B$-modules is
$\tilde{r}_s = r_s \cap \{ 1,\ldots , \tilde{e} \}$. 

For $s=a+1,\ldots , a+\tilde{e}$, the center of $B\tilde{\epsilon }_s$ is
isomorphic to $R$ and $B\tilde{\epsilon }_s$ is isomorphic to one of 
the $R$-orders in Theorem \ref{cycple1} (ii).
%In fact the Brauer tree $\tilde{T}$ of $B$ is a $m$-fold cover of the
%Brauer tree of ${\cal B}$ ramified at the exceptional vertex.
%As above we can partition it into $\tilde{T}_{odd}$ and $\tilde{T}_{even}$.

Let $n':=|\tilde{r}_1| = \frac{|r_1|}{m}$.
For $s=1,\ldots , a$
let $D_s$ be a central $\tilde{Z}_s$-division algebra of index $m$ and 
$\Omega _s$ be its maximal order with prime element $\wp _s$.
Then
$$B \tilde{\epsilon }_s \cong \Lambda (\Omega _s, (d_i, d_{\delta (i) } ,\ldots , d_{\delta ^{n'-1}(i)} ), H_{n'} ). $$ 

Then the Hasse invariant of $D _s$ 
(as defined in \cite[(31.7)]{Rei}) is independent of $s$ and can
be read off from the planar embedded
 Brauer tree together with the Galois action of  $\Gal (k/\tilde{k})
\cong \Gal (K/\tilde{K}) $ on the modular constituents of any exceptional
character in ${\cal B}$:

\begin{theorem}\label{hasse}
Let $\psi $ 
be a $p$-modular constiuent of any of the exceptional characters in ${\cal B}$.
Let $F$ denote the Frobenius automorphism of $k/\tilde{k}$.
Then there is some $r \in \Z $ prime to $m$ 
such that $$\delta ^{n'}  (\psi ) =  F^r(\psi ) \mbox{ where $n' = |\tilde{r}_1| = \frac{|r_1|}{m}$} .$$
Let $t = r^{-1} \in \Z/m\Z $. 
Then for all $s\in \{1,\ldots , a\} $ the Hasse invariant of $D_s$ is $\frac{t}{m}$.
\end{theorem}

\proof
To simplify notation we again assume that all the character degrees $d_i$ are
equal to 1. Then for 
$s\in \{1,\ldots , a\}$ the order 
$B\epsilon _s = \Lambda (\Omega _s,(1,\ldots , 1), H_{n'} ) $ 
and  
$$P = \left( \begin{array}{cccc}
0 &  \ldots & 0 & \wp _s \\
1 & 0 & \ldots &  0 \\
0 & \ddots &  \ddots & \vdots \\
0 & \ldots &  1 & 0 
\end{array} \right)
$$ is a generator of 
$J(B\tilde{\epsilon }_s) $.
Then $P $ also generates the Jacobson radical of 
 ${\cal B} \epsilon _s = R \otimes _{\tilde{R}} B \epsilon _s$.
Let $L_{\psi }$ be a ${\cal B}\epsilon_s$-lattice whose  head has character $\psi $.
Then the head of $L_{\psi } P^{n'} $ has character  
$\delta ^{n'} (\psi )$ which is Galois conjugate to $\psi $ and hence 
of the form $F^r(\psi )$ for some $r$. Therefore conjugation by
$P^{n'} = \diag (\wp _s , \ldots ,\wp _s)$ induces the Galois automorphism 
$F^r$ on the inertia subfield $K$ of $D _s$. 
By the general theory of division algebras over local fields (see \cite{Rei})
 $r$ is prime to $m$ and the Hasse invariant of 
$D_s$ is $\frac{t}{m}$ as stated in the theorem.
\eb

The amalgamations in $B$ are described as in Theorem \ref{cycple}
(see \cite[Theorem 8.10]{Ple 83}),
where now the epimorphisms in (i) are only mappings between $\tilde{R}$-orders.
For (iii) one should note that 
$R/pR \cong \Omega _s/\wp _s \Omega _s $ for all $s=1,\ldots ,a$.

Similarly as in Theorem \ref{unique} one shows:

\begin{theorem}\label{unique2}
The description above (see \cite[p. 140ff]{Ple 83}) determines 
$B$ up to isomorphism.
\\
More precisely let 
$\tilde{\Gamma }_a := (\tilde{\epsilon }_1 + \ldots + \tilde{\epsilon}_a) B$.
and let $\tilde{e}_1,\ldots ,\tilde{e}_{n'}$ $({n'}=|\tilde{r}_1| = \frac{|r_1|}{m} )$
be lifts of the central primitive idempotents in $\tilde{\Gamma }_a/J(\tilde{\Gamma }_a ) $.
Again we assume that the $k$-dimensions of the simple $\tilde{\Gamma }_a$-modules
are 1. Then
$\tilde{e}_i \tilde{\Gamma } _a \tilde{e} _i $ is generated as a
$Z(\tilde{\Gamma} _a)$-order by 
$(\zeta _1,\ldots , \zeta _a)$ and $(\wp _1,\ldots , \wp _a )$, where 
$\zeta _s \in \Omega _s $ is a primitive $(q^m-1)$st root of unity $(q:=|\tilde{k}| =|\tilde{R}/p\tilde{R} |)$
and the prime elements $\wp _s\in \Omega _s$ are chosen such that 
$\zeta _s ^{\wp _s} = F^r(\zeta _s ) = \zeta _s ^{q^r}$ where  
$r$ is as in Theorem \ref{hasse} (i.e. $\frac{t}{m}$ is the Hasse invariant of 
$D_s$ where $rt \equiv 1 \pmod{m} $).
The remaining generators of $\tilde{\Gamma } _a$ are
$\tilde{e}_{i+1,i} \in \tilde{e}_{i+1} \tilde{\Gamma } _a \tilde{e}_{i}$ ($i=1,\ldots , {n'}-1$)
and $\tilde{g}_{1,{n'}} \in \tilde{e}_1 \tilde{\Gamma } _a \tilde{e}_{n'}$ defined analogously
 to the ones 
in Theorem \ref{unique}.
\end{theorem}

\proof
Let $L _i = \tilde{e}_i \tilde{\Gamma }_a$ 
 be any projective indecomposable $\tilde{\Gamma }_a$-lattice 
($i=1,\ldots , n'$).
Then, by the above, the endomorphism ring of
 $L_i$ is a successive amalgam of the 
orders $\Omega _s$, $s=1,\ldots , a$.
Since $\tilde{k}$-automorphisms of $\Omega _s /\wp _s \Omega _s $ lift 
to (inner) $\tilde{R}$-automorphisms of $\Omega _s$, this ring is uniquely determined
by \cite[Theorem 8.10]{Ple 83}  up to isomorphism and 
$$\tilde{e}_i \tilde{\Gamma }_a \tilde{e}_i = \End _{\tilde{\Gamma }_a} (L _i ) 
= \langle (\zeta _1,\ldots , \zeta _a), 
(\wp _1,\ldots , \wp _a ) , Z(\tilde{\Gamma}_a) \rangle .$$
To generate $\tilde{\Gamma }_a$, by Nakayama's lemma,
 it is enough to choose additional 
elements of $\tilde{e}_i \tilde{\Gamma }_a \tilde{e}_j $ ($i\neq j \in \{ 1,\ldots n' \}$) that generate 
$$\tilde{e}_i \tilde{\Gamma }_a \tilde{e}_j /
( \tilde{e}_i J(\tilde{\Gamma }_a)^2 \tilde{e}_j  + p \tilde{e}_i \tilde{\Gamma }_a \tilde{e}_j 
)$$ as an 
$\tilde{e}_i \tilde{\Gamma }_a \tilde{e}_i  $-module.
The same arguments as in the proof of Theorem \ref{unique} now imply 
the theorem.
\eb

%\begin{rem}
%In the analogous notation of Theorem \ref{unique} let 
%$$ 
%(\tilde{x}_1\wp_1 \tilde{e}_{1,{n'}}^1 ,\ldots , \tilde{x}_a\wp_a \tilde{e}_{1,{n'}}^a) =: \tilde{g}_{1,{n'}} \in
%\tilde{e}_1 \tilde{\Gamma }_a \tilde{e}_{n'}$$
%for some units $\tilde{x}_s \in \Omega _s^*$.
%Replacing $\tilde{g} _{1,n'}$ by some other generator, we may assume that
%$\tilde{x}_1 = 1$.
%Then the ${n'}$-th power of $\tilde{e}_{1,2}  + \ldots + \tilde{e}_{{n'}-1,{n'}} + \tilde{g} _{{n'},1} $ is 
%$\wp:=(\wp_1I_{n'} ,\ldots , \tilde{x}_a\wp_a I_{n'})\in \tilde{\Gamma _a} $.
%Since $\wp $ has to induce the same Galois automorphism on $R \subset \Omega _s $, and therefore $\tilde{x}_s $ can be chosen to ly in $\tilde{Z}_s^*$ 
%and even to be in $\tilde{R}^*$.
%\end{rem}

\begin{kor}
The planar embedded Brauer tree 
together with the character fields $K$, 
$\tilde{Z}_1$, $\ldots $, $\tilde{Z}_a $ and the Galois action on the 
modular constituents of the exceptional characters determine the 
block $B$ of $\Z_pG$ up to Morita equivalence.
\end{kor}

\section{The radical idealizer chain for blocks with cyclic defect groups}
\label{head}

In this section we will investigate the radical idealizer chain for blocks
 with cyclic defect group, where we mainly concentrate on describing the head order.
%Starting with an order $\Lambda _0 := \Lambda $ one defines successively
%bigger orders 
%$$\Lambda _0 \subseteq \Lambda _1 \subseteq \ldots \subseteq \Lambda _N = \Lambda _{N+1} $$ defined as 
%$\Lambda _i = \Id (J( \Lambda _{i-1})) $ (see \cite{BeZ})
%the idealizer of the Jacobson radical of $\Lambda _{i-1}$ ($i=1,\ldots $).
%This finite chain will end in a hereditary order $\Lambda _N$ called the
%{\em head order} of $\Lambda $.
Head orders are hereditary orders and hence 
they are the  maximal elements for the ``radically covering''
relation, where an order $\Gamma $ {\em radically covers} and order $\Lambda $,
$\Gamma \succ \Lambda $, 
if $\Gamma \supseteq \Lambda $ and 
 $J(\Gamma ) \supseteq J(\Lambda )$.
Then for all orders in the idealizer chain 
$\Lambda _i \succ \Lambda _{i-1}$
(see \cite[Section 39]{Rei}).
Moreover it is easy to see that if 
$\Gamma \succ \Lambda $ then every simple $\Gamma $ module is 
semi-simple as a $\Lambda $-module (see \cite[Lemma 2.2]{radid}).
In particular the simple $\Lambda _N$-modules are 
semi-simple $\Lambda $-modules.

We will use the notation introduced in the last section and perform the
calculations for the block ${\cal B}$ of $RG$.
The block $B$ of $\Z_pG$ can be treated similarly and the 
head order of $B$ is easily derived from the one of ${\cal B}$ 
(see Remark \ref{headzp}).
However, it is crucial for the whole process that $R$ is an
unramified extension of $\Z _p$.

For the radical idealizer process we treat the exceptional vertex 
$\Gamma _a$ and $\Gamma _0$ separately always keeping track of the 
amalgamations between them, which are controlled by the following lemma.

\begin{lemma}\label{amalgam}
Let $S$ be a discrete valuation ring with prime element $\pi $ and 
let $\Lambda _i$ ($i=1,2$) be $S$-orders.
Given epimorphisms $\varphi _i : \Lambda _i \to X := S^{s\times s} / \pi ^t S^{s\times s} $
let
$$\Lambda := \{ (x_1,x_2)\in \Lambda _1 \oplus \Lambda _2 \mid \varphi _1 (x_1) = \varphi _2 (x_2) \} .$$
Then 
$$\Id (J(\Lambda )) \supseteq 
\{ (x_1,x_2)\in \Lambda _1 \oplus \Lambda _2 \mid \overline{\varphi _1 (x_1) }= \overline{\varphi _2 (x_2)} \}  =: \Gamma $$
where $\overline{\phantom{s}} : X \to 
S^{s\times s} / \pi ^{t-1} S^{s\times s} $ is the natural epimorphism.
\end{lemma}

\proof
Clearly $J(\Lambda ) = 
 \{ (y_1,y_2)\in J(\Lambda _1) \oplus J(\Lambda _2) \mid \varphi _1 (y_1) = \varphi _2 (y_2) \} $ and $\varphi _i (J(\Lambda _i)) = J(X) = \pi X $ for $i=1,2$.
Let $(x_1,x_2)\in \Gamma $ and $(y_1,y_2) \in J(\Lambda )$.
Then clearly $x_i y_i$ and $y_i x_i$ are in $J(\Lambda _i)$ ($i=1,2$).
Since $\varphi _1$ is surjective, there is $z_1\in \Lambda_1$ with 
$\pi \varphi _1(z_1) = \varphi _1 (y_1)$.
Choose $z_2\in \Lambda _2$ with $\varphi _2(z_2) = \varphi_1(z_1)$.
Then 
$$\varphi _1(y_1x_1) =  \varphi _1(z_1) \pi \varphi _1(x_1)
= \varphi_2 (z_2) \pi \varphi _2(x_2) = \varphi _2 (y_2x_2)$$
and similarly $\varphi _1 (x_1 y_1) = \varphi _2(x_2y_2 )$.
Hence $(x_1,x_2) \in \Id(J(\Lambda ))$.
\eb

The following trivial lemma suffices to deduce the 
head order of $\Gamma _a$.

\begin{lemma}
Let $\Lambda $ be an order in ${\cal A}$ and $\epsilon $ a 
central idempotent of ${\cal A} $.
Then $$\Lambda \epsilon \subseteq \Id (J(\Lambda )) \epsilon \subseteq
\Id (J(\Lambda \epsilon )) .$$
\end{lemma}

\begin{kor}\label{ausnahme}
The head order of $\Gamma _a$ is 
$\oplus_{s=1}^a {\cal B} \epsilon _s $.
\\
Similarly the 
head order of $\tilde{\Gamma }_a$ is 
$\oplus_{s=1}^a B \tilde{\epsilon }_s $.
\end{kor}

\proof
The orders 
${\cal B} \epsilon _s = \Id (J({\cal B} \epsilon _s))$ and
$B \tilde{\epsilon }_s $ are already
hereditary for $s=1,\ldots , a$.
\eb

Note that this corollary is not true, when $R$ is replaced by a
ramified extension of $\Z_p$.

\subsection{The first steps.}

The main task to calculate the idealizer chain for $\Gamma _0$
is to calculate the one of $\epsilon _s \Gamma _0$ for 
$s=a+1,\ldots , a+e$. These orders have a certain symmetry with respect to
a cyclic permutation of their simple modules, and therefore can be 
encoded in a simple way.
All orders in this radical idealizer chain share this symmetry.

\begin{definition}
For $v= (v_0,\ldots , v_{n-1} ) \in \Z ^n $ and 
$d=(d_1,\ldots , d_n) \in \N ^n$ define
$$\Lambda (d,v) := \Lambda (R,d,M) := \{ X = (x_{ij}) \in K^{D\times D } \mid x_{ij} \in \pi ^{m_{ij}} R^{d_i\times d_j}  \} $$ where
$$m_{ij} = \left\{ \begin{array}{ll}
v_{j-i} &  \mbox{ if } j\geq i  \\ 
v_{n+j-i}-v_{n-1} &  \mbox{ if } j < i \end{array} \right. $$
and $D=\sum _{i=1}^n d_i$.
\end{definition}

\begin{rem}
Since the dimension vector $d$ will be fixed most of the time, 
we will omit it and let $\Lambda (v_0,\ldots , v_{n-1}) := \Lambda (d,v)$.
\end{rem}

The order $\Gamma _0$ is an amalgam of the orders 
$\epsilon _s {\cal B} = \epsilon _s \Gamma _0 $ $(s=a+1,\ldots , a+e )$
of the form 
$ \Lambda (R,d, a H_n) \cong \Lambda (0,\underbrace{a,\ldots ,a}_{n-1} )
$ for some dimension vector $d$ and $n=|r_s|$. 
The amalgamations in $\Gamma _0$ are only on the diagonal, more 
precisely, the part of ${\cal B} $ belonging to $\epsilon _s {\cal B} $ 
is of the form 
$$\Lambda (R,d , 
\left( \begin{array}{ccccc}
\underline{0}_a & a & \ldots & \ldots & a \\
0 & \underline{0}_a & a & \ldots & a  \\
\vdots & \ddots & \ddots & \ddots & \vdots \\
0 & \ldots & \ldots & \underline{0}_a & a \\
0 & \ldots & \ldots & 0 & \underline{0}_a \end{array} \right) 
) = \Lambda (\underline{0}_a, a ^{n-1}) $$
where all the underlined entries obey a certain congruence modulo $p^a$ to
a diagonal entry in some other $\epsilon _t {\cal B} $ ($t\neq s $) which
is indicated by underlining the $0$ and the index $a$.
By Lemma \ref{amalgam} these amalgamations will decrease by 1 in each step until
after $a$ steps the order ${\cal B}_a$ contains  the central primitive
idempotents $\epsilon _{a+1},\ldots , \epsilon _{a+e}$.

In the following we fix some
 $s\in \{a+1,\ldots , a+e\} $, put $n:= | r_s | $, and let 
$$\Lambda := \Lambda _0 := {\cal B}\epsilon _s \subseteq
\Lambda _1 := {\cal B}_1 \epsilon _s \subseteq \ldots  \subseteq 
\Lambda _N := {\cal B}_N \epsilon _s  $$ 
where $${\cal B}=: {\cal B}_0 \subset {\cal B}_1 \subset \ldots \subset {\cal B}_N  = {\cal B}_{N+1}$$
is the radical idealizer chain of ${\cal B}$.
Together with the structure of $\Lambda _i$ we keep track of 
the additional information, how ${\cal B}_i$ is embedded into the direct sum of
the $ {\cal B}_i \epsilon _s$ using the notation above.

\begin{lemma} 
If $a\geq n$ then
$$\Lambda _{n}  \cong \Lambda (\underline{0}_{a-n} , (a-n)^{n-1} ).$$
\end{lemma}

\proof
An easy induction on $j$  shows that for $j=1,\ldots , n$ 
$$\Lambda _{j}  = 
\Lambda (\underline{0} _{a-j} , (a-j+1),(a-j+2),\ldots , a-1,a^{n-j} ) .$$
Then 
$$\Lambda _n  
=\Lambda (\underline{0} _{a-n} , (a-n+1),(a-n+2),\ldots , a-1 ) 
\cong  \Lambda (\underline{0}_{a-n} , (a-n)^{n-1} ) $$
by conjugation with the diagonal matrix 
$\diag(1,\pi , \pi^2 ,\ldots , \pi ^{n-1} ) .$
\eb

Inductively we get 

\begin{kor}\label{redb}
Let $a = z_s n + b $ with $0\leq b < n$ and $m_{0} := z_s n+1$.
Then $$\Lambda_{m_{0} -1 }  = \Lambda (\underline{0}_b, b^{n-1} ). $$
\end{kor}

If $b = 0$ then $\Lambda  _{m_{0}-1}$ is already a maximal order and we are done.

\begin{lemma}\label{anfang}
%Let $\Lambda _{m_{0}}  = \Lambda ((\underline{0}_b ,b^{n-1})) $ 
%with $0<b<n$ and
%$\Lambda_{m+m_{0}} := {\cal B}_{m+m_0} \epsilon _s  $ $m=1,2,\ldots $ as above.
Assume that $b>0$ and define $l_0,x_0$ by  $n=l_0b +x_0 $ with $0 < x_0 \leq b$.
Then 
$$\Lambda _{m_{0}}  = \Lambda (\underline{0}_{b-1} ,b^{n-1}) .$$
If $0\leq m<n-l_0$ then
\begin{itemize}
\item[(a)] 
$\Lambda _{m+m_{0}} $ is of the form $\Lambda _{m+m_{0}} = \Lambda (\underline{0}_{f(m)},v_1,\ldots , v_{n-1} )$
with $0 < v_1 \leq v_2 \leq \ldots \leq v_{n-1} = b $, where
$f(m) = \max\{0,b-m-1\}$.
\item[(b)]
If $m=l(b-1) +y $ with $0\leq y < b-1$ then 
$$\Lambda _{m+m_{0}} = \Lambda (
\underline{0}_{f(m)} , 1^l,2^l,\ldots , (b-y-1)^l,(b-y)^{l+1},\ldots ,  (b-1)^{l+1},b^{n-m-1})$$
i.e. 
$\Lambda_{m+m_{0}} = \Lambda (\underline{0}_{f(m)},v_1,\ldots , v_{n-1})$ with 
$$ v_j = \left\{ \begin{array}{ll} \lfloor \frac{j-1}{l} \rfloor +1 & \mbox{ if }
1\leq j \leq (b-y-1)l \\
b-y+\lfloor \frac{j-1-(b-y-1)l}{l+1} \rfloor & \mbox{ if } (b-y-1)l < j \leq (b-1)l+y \\
b & \mbox{ if } j > (b-1) l + y \end{array} \right. $$
\item[(c)] 
The radical $J(\Lambda _{m+m_{0}}) = \Lambda (\underline{1}_{f(m)},v_1,\ldots , v_{n-1})$.
\end{itemize}
\end{lemma}

\proof
The form of $\Lambda_{m_0}$ is clear.
For the other statements
we argue by induction on $m$, where the case $m=0$ is trivial.
Assume that $m<n-l_0-1$ and that $\Lambda _{m+m_{0}} =\Lambda (\underline{0}_{f(m)},v_1,\ldots , v_{n-1})$
 has the properties $(a),(b),(c)$.
Then $\Lambda _{{m+m_{0}}+1} $ is of the
form $$\Lambda _{{m+m_{0}}+1} = \Lambda (\underline{0}_{f(m+1)},\tilde{v}_1,\ldots \tilde{v} _{n-1}) ,$$
since the inequalities on the entries of the exponent matrix preserve the
symmetry conditions in $(a)$.
The form of the amalgamations follows from Lemma \ref{amalgam}.
Clearly $\tilde{v}_i \leq v_i$ for all $i$ and  $\tilde{v}_0 = 0$.
The remaining conditions in (a) and the property (c) follow once we have shown 
(b). 
Let $v_i':= v_i$ for $i>0$ and $v_0':= 1 = v_0+1$.
Then the conditions on $m_{1j} = \tilde{v}_{j-1}$ ($j>1$)
that $\Lambda _{m+m_0+1}$ lies in the left idealizer of 
$J(\Lambda _{m+m_{0}})$ read as 
$$\tilde{v}_{j-1} \geq \max \{ v_{k-1} - v_{k-j}' \mid k=j,\ldots , n \}  =: \mbox{ $\max _1$} $$ and
$$\tilde{v}_{j-1} \geq \max \{ b+v_{k-1}' - v_{k+n-j} \mid k=1,\ldots , j-1 \}  =: \mbox{$ \max _2$} .$$
The inequalities for the right idealizer of $J(\Lambda _{m+m_{0}})$ read as
$$\tilde{v}_{j-1} \geq \max \{ 
v_{j-1} -1, v_{j-k_1}' - v_{n+1-k_1} + b,
v_{n+j-k_2}-v_{n+1-k_2} \mid
 k_1=2,\ldots , j , 
 k_2=j+1,\ldots , n \} $$
and agree with the conditions above after an easy variable transformation.
Hence right and left idealizer of 
$J(\Lambda _{m+m_{0}})$ coincide and are equal to $\Id (J(\Lambda _{{m+m_{0}}}))$.

By the induction assumption for all $1\leq i \leq n$
$$ v_{i-1} = \left\{ \begin{array}{ll} \lfloor \frac{i-2}{l} \rfloor +1 & \mbox{ if }
1\leq i-1 \leq (b-y-1)l \\
b-y+\lfloor \frac{i-2-(b-y-1)l}{l+1} \rfloor & \mbox{ if } (b-y-1)l < i-1 \leq (b-1)l+y \\
b & \mbox{ if } i-1 > (b-1) l + y . \end{array} \right. $$

Since the `slope' of $v$ is decreasing
$v_{k-1} - v_{k-j}'$ is maximal if $v_{k-j}'$ is the last 1 in $v'$,
hence if $k-j = l$ i.e. $k=l+j$. 
If $k:=\min (l+j,n)$ then 
$$\mbox{$ \max _1 $} = v_{k-1} - 1 = \left\{ \begin{array}{ll} 
\lfloor \frac{j-2}{l} \rfloor + 1 & \mbox{ if } 1<j-1 \leq (b-y-2)l    \\
b-y + \lfloor \frac{j-2-(b-y-2)l}{l+1} \rfloor - 1 & \mbox{ if }
(b-y-2)l <  j-1 \leq (b-2)l+y   \\
b-1 & \mbox{ if }  (b-2)l+y < j-1  .
\end{array} \right. $$
This implies that 
$\max _1 = v_{j-1}$ if $j-1 \leq (b-y-2)l$.
If $(b-y-2)l < j-1 \leq (b-y-1)l$ then 
$$v_{j-1} = \lfloor \frac{j-2}{l} \rfloor +1 = b-y-1 = \mbox{$ \max_1  $}
= b-y + \lfloor \frac{j-2-(b-y-2)l}{l+1} \rfloor -1 .$$

If $(b-y-1)l < j-1 \leq (b-2)l + y $ then 
$\max _1 = b-y-1 + \lfloor \frac{j-2-(b-y-2)l}{l+1} \rfloor $
and $v_{j-1} = b-y + \lfloor \frac{j-2-(b-y-1)l}{l+1} \rfloor $.
Therefore $\max _1 < v_{j-1} $ if and only if $j-3-y+b$ is divisible by
$l+1$, i.e. 
$$j-1 = (b-2)l+y-x(l+1), \ x=0,1,\ldots , y-2 $$
when $v_{j-1} = b-x-1$ is the first occurrence of $b-x-1 $ in $v$.

If $(b-2)l +y < j-1 \leq (b-1)l+y $ then 
$v_{j-1} = b-1 = \max _1 $ and if $j-1 > (b-1)l+y$ then $\max _1 = b-1 <
v_{j-1} = b$. 

For $\max _2$ one finds that $b+ v_{k-1}' -v_{k+(n-j)}$ is maximal
if $k=j-1$ since the `slope' of $v$ is decreasing.
Hence $$\mbox{$ \max _2 $} = b+v_{j-2}' - v_{n-1}  = v_{j-2} ' .$$

Combining these conditions % coming from $\max _1$ and $\max _2$ 
one finds that 
$$\tilde{v}_{j-1} = \left\{ \begin{array}{ll}
v_{j-1} - 1 & \mbox{ if } j-1 = (b-2) l + y - x(l+1) \mbox{ for } x=0,\ldots , y-1 \\
v_{j-1} & \mbox{ otherwise } \end{array} 
\right. .$$

With these $\tilde{v}_j$ the multiplication by $\Lambda _{m+m_0+1} = 
\Lambda (\underline{0}_{f(m+1)},\tilde{v}_1,\ldots ,\tilde{v}_{n-1})$ 
preserves the congruences in $J(\Lambda _{m+m_0})$ given in (c), since
$$b-m-1\leq f(m)  \leq \tilde{v}_{j-1} + v_{n+1-j} - b 
\mbox{ and  } b-m-1 \leq f(m) \leq \tilde{v}_{n+1-j} -b + v_{j-1} $$ 
for all $j$.
This implies part (b) of the lemma.
\eb

\begin{kor}\label{defm1}
Let $m_1:=n-l_0-1+m_{0}$ and $y = x_0-1$. Then
 $$\Lambda _{m_1} = 
\Lambda (0,1^{l_0},\ldots , (b-y-1)^{l_0},(b-y)^{{l_0}+1},\ldots , (b-1)^{{l_0}+1},b^{l_0})   =: \Lambda (v^{(1)}).$$
\end{kor}

\begin{kor}
For all $s\in \{ 1,\ldots , a+e \}$ 
the $s$-th component of the head order of ${\cal B}$ is
equal to the head order of the projection $\epsilon _s {\cal B}$.
\end{kor}

%Another corollary from the proof of Lemma \ref{anfang} is the following 
%observation.

%\begin{kor}
%For blocks ${\cal B}$ with cyclic defect groups, 
%left and right idealizer chains coincide and
%hence are equal to the idealizer chain of ${\cal B}$.
%\end{kor}

\subsection{The head order.}

If $n=(l_0+1)b$ is divisible by $b$,
then the order $\Lambda _{m_1}$ as defined in Corollary \ref{defm1}
is already hereditary.
More precisely 
we have the following 

\begin{lemma}{\label{ggt=b}}
Let $b$ be a factor of $n = lb$. 
Then $\Lambda _{m_1} = \Lambda (0,1^{l},\ldots , (b-1)^{l} , b^{l-1}) $ 
is  hereditary, 
$\Lambda _{m_1} \sim \Lambda (R,(D_1,\ldots ,D_{l}), H_{l})$
where $D_i = \sum _{j = 0}^{b-1} d_{j l + i} $.
\end{lemma}

\proof
$\Lambda _{m_1} = \Lambda (R,d,M)$ where 
$$m_{ij} = \lfloor \frac{j-i-1 }{l} \rfloor + 1 .$$
Let $t_i :=  m_{i1} =  \lfloor \frac{-i}{l} \rfloor - 1 $.
Conjugating by the diagonal matrix $T:=\diag (\pi ^{t_i})$ one obtain the
conjugate order 
$\Lambda _{m_1} ^{T} = \Lambda (R,d, \tilde{M} ) $, where 
$$\tilde{m}_{ij} = m_{ij} - t_i + t_j = 
\lfloor \frac{j-i-1 }{l} \rfloor + 1 -  \lfloor \frac{-i}{l} \rfloor 
+  \lfloor \frac{-j}{l} \rfloor  .$$ 
Writing $j = j_1 l + j_2$ and $i  = i_1 l + i_2 $ with $0 < j_2,i_2 \leq l$
 one gets 
$$\tilde{m}_{ij} 
\lfloor \frac{j_2-i_2-1 }{l} \rfloor + 1 -  \lfloor \frac{-i_2}{l} \rfloor 
+  \lfloor \frac{-j_2}{l} \rfloor  
=\lfloor \frac{j_2-i_2-1 }{l} \rfloor + 1  = \left\{ \begin{array}{ll} 
0 & \mbox{ if } j_2 \leq i_2 \\
1 & \mbox{ if } j_2 > i_2 \end{array} \right.
.$$
Hence after reordering the constituents $\Lambda _{m_1}$ has the form as
claimed in the lemma.
\eb

We now assume that $0<x_0<b$. 
Continuing to
trace down the radical idealizer process like in Lemma \ref{anfang}
  seems to be a rather tedious work.
If $x_0\geq \frac{b}{2}$ then after $m_2 = b-x_0-1$ steps one arrives at an
order 
$$\Lambda _{m_1+m_2} = 
\Lambda (0,1^{l_0} , 2^{l_0+1}, 3^{l_0},4^{l_0+1} ,\ldots , z^{l_0} ,
(z+1)^{l_0+1}, (z+2)^{l_0+1}, \ldots , (b-1)^{l_0+1}, b^{l_0})$$
where $z= 2b-2x_0-1 $.
If $x_0 \leq \frac{b}{2}$ then after $m_2 = x_0-1$ steps one arrives at an
order 
$$\Lambda _{m_1+m_2} = 
\Lambda (0,1^{l_0} , 2^{l_0}, \ldots , z^{l_0} ,
(z+1)^{l_0+1}, (z+2)^{l_0}, \ldots , (b-2)^{l_0+1}, (b-1)^{l_0}, b^{l_0})$$
where $z= b-2x_0+1 $.
If $x= \frac{b}{2} = \gcd(n,b) = d $ then $\Lambda _{m_1+m_2}$ is again
hereditary 
$$\Lambda _{m_1+m_2} \sim \Lambda (R,(D_1,\ldots ,D_{2l_0+1}), H_{2l_0+1})$$
where $D_j = \sum _{i \equiv -l_0 j} d_i  $, where the
congruence is modulo $2l_0+1 = \frac{n}{d} $.

Instead of continuing like this, 
we prefer to calculate the head order $\Lambda _N = \epsilon _s {\cal B}_N $, 
which is also the head order of $\Lambda _{m_1}$,
directly where we need the following trivial lemma:

\begin{lemma}{\label{cov}}
Let $\Lambda \subset \Gamma $ be two orders with $J(\Lambda ) \subset J(\Gamma )$.
If $e\in \Lambda $ is an idempotent then $J(e\Lambda e) \subset J(e\Gamma e)$.
\end{lemma}

\proof
$$J(e\Lambda e) = e J(\Lambda ) e  \subseteq e J(\Gamma ) e = J(e \Gamma e ). $$
\eb

The  head order $\Lambda _N$ of $\Lambda _{m_1}$ has the following properties:

\begin{properties}\label{prop}
{\rm 
\begin{itemize}
\item[0)]
$\Lambda _N $ is of the form $\Lambda (w)$ for some $w\in \Z_{\geq 0}^n$.
\item[1)] $\Lambda _N$ is an order, i.e. for all $i<j<k$ one has
\begin{itemize}
\item[(i)] $b-w_{n+j-k} \leq w_{k-i} - w_{j-i} \leq w_{k-j}  $
\item[(ii)] $b-w_{n+i-k} \leq w_{k-j} - w_{n+i-j} + b  \leq w_{k-i}  $
\item[(iii)] $b-w_{n+i-j} \leq w_{n+j-k} - w_{n+i-k} \leq w_{j-i}  $
\end{itemize}
which just expresses the fact that the entries $m_{ij}$ in the exponent matrix 
of $\Lambda _N$ satisfy
$m_{ik} + m_{kj} \geq m_{ij} $ for all $i,j,k \in \{ 1,\ldots , n\}$.
\item[2)] $\Lambda _N$ is hereditary, i.e. 
$$w_{j-1} + w_{n+1-j} - b  \in \{ 0,1 \} 
\mbox{ for all } j > 1 .$$ 
\item[3)] $\Lambda _N$ radically covers the order
$\Lambda _{m_1} = \Lambda (v^{(1)})$ defined in  Corollary \ref{defm1}.
 %which means that 
%$\Lambda _N \supseteq \Lambda _{m_1}$ and 
%$J(\Lambda _N )\supseteq J(\Lambda _{m_1})$.
This property implies with Lemma \ref{cov} 
that $w_{j-1} = v^{(1)}_{j-1}$ and $w_{n-j+1} = v^{(1)}_{n-j+1}$ whenever
$v^{(1)}_{j-1} + v^{(1)}_{n+1-j} - b  = 1 $. 
In particular 
\item[3')] 
$w_1=\ldots = w_{l_0} = 1 $, 
$w_{n-1}=\ldots = w_{n-l_0} = b $, and 
$ w_{n-l_0-1} = b-1 $.
\end{itemize}
}
\end{properties}

\begin{lemma}\label{bteiltnichtn}
$\Lambda _N $ is uniquely determined by Properties \ref{prop} 
0), 1), 2), and 3').
More precisely
let  $n=l_0b+x_0$ be as above and assume that $1\leq x_0 \leq b-1$. Then
\begin{itemize}
\item[(i)]
$\Lambda _N = \Lambda (w)$ where $
w = (0,1^{l_1} , 2^{l_2} \ldots , b^{l_b} ) $ with
$l_1=l_0=l_b$ and $l_j \in \{l_0,l_0+1\}$ for all $j=1,\ldots , b$.
\item[(ii)]
Let $e:=(e_1,\ldots , e_b)$, where $e_k = l_k-l_0 \in \{ 0,1 \}$ 
for $k=1,\ldots b-1$ and $e_b:= 1$.
For all $j$ let $a_j := \sum _{k=1}^j e_k $.
Let $d:= \gcd (n,b) = \frac{b}{i}$.
Then $x_0 = \frac{a_i}{i} b$, $e_i = 1$ and $$e = (e_1,\ldots , e_i)^d = 
(e_1,\ldots , e_i,e_1 ,\ldots , e_i , \ldots , e_1,\ldots , e_i ).$$
The entries of $w$ are uniquely determined by 
$$a_j = \lfloor \frac{x_0\cdot j}{b} \rfloor \mbox{ for all } j=1,\ldots , b .$$
\end{itemize}
\end{lemma}

\proof
(i)
Property \ref{prop} 3') 
together with Property \ref{prop} 1) (i) 
(for $i=1$)  show that for $0 < k-j \leq l_0$ 
$$ 0 \leq w_{k-1} - w_{j-1} \leq 1 $$ and if $k-j \geq l_0+1$, then 
$ w_{k-1} - w_{j-1} \geq 1 $.
This implies (i).
\\
(ii)
Put $d = \gcd(b,x_0) = \frac{b}{i}$.
Then $ i = \min \{ j \in \{ 1,\ldots ,b \} \mid \frac{b}{j} \mbox{ divides } x_0 \}$.
\\
We now show by induction on $j$ that 
$a_j = \lfloor \frac{x_0\cdot j}{b} \rfloor $ and $l_j = l_{b-j+1}$
for $j=1,\ldots , i-1 $.
This is clear for $j=1$ since $l_1=l_0=l_b$ 
and $a_1 = 0 = \lfloor \frac{x_0}{b} \rfloor $.
Assume that $1< j\leq  i-1 $ and that 
$a_{k} = \lfloor \frac{x_0\cdot {k}}{b} \rfloor $ and 
$l_{k} = l_{b-k+1}$
for $k=1,\ldots , j-1 $.
Let $$X_1:=(  (t-1),t^{l_t} , \ldots , (t+j-1)^{l_{t+j-1}} ,(t+j)  ) $$ 
be a {\em subsequence} of $w$.
Then the difference between the first and the last entry of $X_1$ is
$j+1$ and the distance between these entries is
$\sum _{q=t}^{t+j-1} l_q+1 $.
Since $w_{l_0j+a_j} = j$, Property \ref{prop} 1) (i) implies that 
$$\sum _{q=t}^{t+j-1} l_q  
= l_0j + \sum _{q=t}^{t+j-1} e_q  
\geq l_0j + a_j  \mbox{ for all } 1\leq t < b+1-j .$$
Similarly for subsequences of $(w,w)$ of the form
$$X_2:=( (b-t-1),(b-t)^{l_{b-t}}, \ldots , b^{l_b} , 0, 1^{l_1} ,\ldots,
(j-t-1)^{l_{j-t-1}} , (j-t)) $$  
Property \ref{prop} 1) (ii) implies that 
$$\sum _{q=b-t}^{b} l_q +1 + \sum_{q=1}^{j-t-1} l_q =
l_0j + \sum _{q=b-t}^{b} e_q + \sum_{q=1}^{j-t-1} e_q 
 \geq l_0j + a_j  \mbox{ for all } 0\leq t \leq j .$$
This implies that for every subsequence of length $j$ of the sequence 
$(e,e)$, the sum over the entries in this subsequence is $\geq a_j$ and
therefore
 $$x_0 = \sum _{t=1}^{b} e_t    \geq  \frac{b}{j} a_j .$$
Let $b_j := \sum _{t=b-j+1}^b e_t $.
Similar arguments as above, using the second and second last entries of the 
sequences $X_1$ and $X_2$ above and the fact that 
$w_{n-jl_0-b_j} = b-j $, show that 
$$\sum _{q=t}^{t+j-1} l_q   
= l_0j + \sum _{q=t}^{t+j-1} e_q    
\leq l_0j + b_j  \mbox{ for all } 1\leq t \leq b+1-j $$
and 
$$\sum _{q=b-t}^{b} l_q +1 + \sum_{q=1}^{j-t-1} l_q =
l_0j + \sum _{q=b-t}^{b} e_q + \sum_{q=1}^{j-t-1} e_q 
 \leq l_0j + b_j  \mbox{ for all } 0\leq t \leq j $$
which yields
$$\frac{b}{j} a_j \leq x_0  \leq \frac{b}{j} b_j  .$$ 

By induction hypothesis, we have 
$b_j = a_j+1$ (if $l_{b-j+1} = l_j $) or $b_j = a_j $ (if $l_{b-j+1} = l_0 $ and
$l_{j} = l_0+1 $). 
Note that the case $l_{b-j+1} = l_0+1$ and $l_j=l_0$  is not possible
since then $w_{n-jl_0-b_j} + w_{jl_0+b_j} = b-j+1 + j+1 = 2 $ contradicting Property 
\ref{prop} 2).
If $b_j = a_j$ then  $x_0 = \frac{b}{j} a_j$ and $\frac{x_0j}{b} $ is an integer 
showing that $j\geq i$.
If $b_j = a_j + 1 $ then $l_{b-j+1} = l_j $ 
and 
$$ \frac{jx_0}{b} -1 \leq a_j \leq \frac{jx_0}{b} $$
which give $a_j = \lfloor \frac{jx_0}{b} \rfloor $ as claimed, since 
$j \leq i-1$ and hence $\frac{jx_0}{b} $ is not an integer.

It remains to show that if $j=i$, i.e. $\frac{b}{j} = \gcd (b,x_0) = \gcd(b,n) $,
then $a_j = a_i = \frac{ix_0}{b}$ and $e$ and $\Lambda _N$ are as claimed.
For this it is enough to show that $a_i = b_i$, since then every 
subsequence of $e$ of length $i$ contains exactly $a_i$ times 1.
Applying this to $(e_1,\ldots , e_i)$ and $(e_2,\ldots , e_{i+1})$ 
this shows that $e_{i+1} = e_1$. Repeating it follows that 
$e=(e_1,\ldots , e_i , e_1,\ldots, e_i, \ldots e_1,\ldots , e_i)$ as claimed. 

Assume that $a_i \neq b_i$. Then $b_i = a_i +1 $ and either 
$a_i = \frac{x_0}{d} $ and $b_i = \frac{x_0}{d} + 1 $ or
$a_i = \frac{x_0}{d} -1$ and $b_i = \frac{x_0}{d}  $ (where $d:=\frac{b}{i} = \gcd(n,b)$). 
Assume the latter, then 
$$x_0 = \sum _{j=1}^b e_j = \sum _{k=0} ^{d-1} \sum _{j=ki+1}^{ki+i} e_j 
\leq d b_i  = x_0 .$$
Hence for all $k$ the sum
$ \sum _{j=ki+1}^{ki+i} e_j  = b_i$, in particular 
$a_i = \sum _{j=1}^i e_j = b_i$. 
In the other case one argues similarly using $a_i$ instead of $b_i$.
\eb

%\begin{prop}\label{descw}
%Let $n=lb+x$ with $0\leq x < b$ 
%(so $l=l_0$, $x=x_0$ if $b$ is not a factor of $n$ and $l=l_0+1$ if $n=lb$).
%Then $\Lambda _N = \Lambda (w)$ with 
%$$w_j = 1 + \lfloor \frac{j-1-\lfloor \frac{xj}{n} \rfloor}{l} \rfloor .$$
%\end{prop}
%
%\proof
%If $x=0$, then $\Lambda (w) = \Lambda _{m_1}$ from Lemma \ref{ggt=b}.
%Assume that $x>0$.
%By Lemma \ref{bteiltnichtn} it is enough to show that 
%$\Lambda (w)$ with $w_j$ as in the proposition 
%is a hereditary order 
%(this will be shown in the proof of the next theorem) 
%that has Property \ref{prop} 3') (which is a straightforward calculation).
%\eb

\begin{theorem}\label{main2}
The head order of ${\cal B} $ is 
$${\cal B}_N = \bigoplus _{s=1}^{a+e} \Delta _s  $$
where $\Delta _s = {\cal B} \epsilon _s 
 \mbox{ for } s=1,\ldots , a $. 
\\
If $s\in \{ a+1,\ldots , a+e \}$ 
let $d:=\gcd (|r_s| , a )$,  $t:= \frac{|r_s|}{d}$ and 
$c:= (\frac{a}{d} )^{-1} \in (\Z/t\Z )^*$.
As in Theorem \ref{cycple1} let
$\sigma := \delta _{|r_s}$ if $s\in T_{even}$ and $\sigma := \rho _{|r_s}$ if $s\in T_{odd}$.
Then the order of $\sigma $ is $|r_s|$ and we define $\tau := \sigma ^t$ and
 $\gamma := \sigma ^c$ and choose $i\in r_s$ arbitrarily.
Then
$$\Delta _s \cong \Lambda (R, (D_i, D_{\gamma  (i) } ,\ldots , D_{\gamma ^{t-1}(i)}) ,  H_{t} ) $$
where $D_j = \sum _{l=0}^{d-1} d_{\tau ^l(j)} $.
\end{theorem}

\proof
For $1\leq s \leq a$ the theorem follows from Corollary \ref{ausnahme}.
For $a+1 \leq s\leq a+e $ let  $n:=|r_s|$, $a=\mu n + b$ with 
$0\leq b<n$. If $b=0$, then $\Delta _s \cong \Lambda _{m_1}$
as defined in Corollary \ref{defm1} is already a maximal order and
the theorem follows from Lemma \ref{ggt=b}.

So assume that $1\leq b \leq n-1$. Then  $d=\gcd (a,n) = \gcd(b,n)$ 
and we write $n=lb+x$ with $0\leq x < b$ and put 
$n=n' d, \ b=b'd, \ x = x' d $.
Then there is $k\in \Z$ with $cb' = 1+n'k$ where $c$ is as defined in the theorem.
For $j\in \Z $ put 
$$f(j): = 1 + \lfloor \frac{j-1-\lfloor \frac{xj}{n} \rfloor}{l} \rfloor 
 = 1 + \lfloor \frac{-1-\lfloor \frac{(x'-n')j}{n'} \rfloor}{l} \rfloor 
.$$
Since $x'-n'=-b'l$ is divisible by $l$, one finds that 
$$f(j+n') = f(j) +b' \mbox{ for all } j\in \Z . $$

Let $$\Lambda := \Lambda (f(0),\ldots , f(n-1)) = 
\Lambda (R,d,M) \mbox{ where }
m_{ij} = f(j-i) .$$

We claim that $\Lambda = \Delta _s$.
By Lemma \ref{bteiltnichtn} it is enough to show that 
$\Lambda $ is a hereditary order that has property \ref{prop} 3').
The latter is checked by a straightforward calculation.
We show that $\Lambda $ is hereditary, by 
establishing an isomorphism with the hereditary order in the theorem.

Put $t_i :=  m_{i1} = f(1-i) $.
Conjugating by the diagonal matrix $T:=\diag (\pi ^{t_i})$ one obtains the
conjugate order
$$\Lambda ^{T} = \Lambda (R,d, \tilde{M} ) \mbox{, where }
\tilde{m}_{ij} = m_{ij} - t_i + t_j  = f(j-i) - f(i) +f(j) .$$
Writing $j = 1+cj_2+n'j_1$ and $i  = 1+ci_2+n'i_1 $ with $0 \leq j_2,i_2 <n'$
 one gets
$\tilde{m}_{ij} = f(c(j_2-i_2)) - f(-ci_2) + f(-cj_2) $ $$ =
 1 + \lfloor \frac{-1-\lfloor \frac{-cb'l(j_2-i_2)}{n'} \rfloor}{l} \rfloor 
  - \lfloor \frac{-1-\lfloor \frac{cb'li_2}{n'} \rfloor}{l} \rfloor 
  + \lfloor \frac{-1-\lfloor \frac{cb'lj_2}{n'} \rfloor}{l} \rfloor  .$$
Since $cb' =1-kn'$ one gets 
$$\tilde{m}_{ij}  =
 1 + \lfloor \frac{-1-\lfloor \frac{l(i_2-j_2)}{n'} \rfloor}{l} \rfloor 
  - \lfloor \frac{-1-\lfloor \frac{li_2}{n'} \rfloor}{l} \rfloor 
  + \lfloor \frac{-1-\lfloor \frac{lj_2}{n'} \rfloor}{l} \rfloor  .$$
Now $0 \leq i_2 < n'$ implies that 
$0\leq \lfloor \frac{li_2}{n'} \rfloor  \leq l-1 $ and therefore 
  $ \lfloor \frac{-1-\lfloor \frac{li_2}{n'} \rfloor}{l} \rfloor  = -1$.
Similarly 
  $ \lfloor \frac{-1-\lfloor \frac{lj_2}{n'} \rfloor}{l} \rfloor  = -1$.
For the first term we have 
$1-n' \leq i_2 - j_2 \leq n'-1 $ implying that 
 $$ \lfloor \frac{-1-\lfloor \frac{l(i_2-j_2)}{n'} \rfloor}{l} \rfloor \in \{ 0,-1 \}.$$
More precisely this yields 
$$\tilde{m_{ij}} = \left\{ \begin{array}{ll} 0 & \mbox{ if } i_2 \geq j_2 
\\ 1  & \mbox{ if } i_2 < j_2  . \end{array} \right.$$
In particular $\Lambda $ is a hereditary order and hence $\Lambda = \Delta _s$.
After a suitable reordering of the constituents the order  $\Lambda ^T \cong 
\Delta _s$ has the form as claimed in the theorem.
\eb

\begin{rem}{\label{headun}}
Let 
$s\in \{ a+1,\ldots , a+e \}$ and $\Delta _s := \epsilon _s {\cal B}_N$.
Let $n:=|r_s| = n' d$,  $d=\gcd (a,n)$, $a=a'd$, and $c a' \equiv 1 \pmod{n'}$.
Let $\nu : \Z/n\Z \to \Z/ n'\Z $ be the natural epimorphism.
Assume that the simple ${\cal B}\epsilon_s$-modules are labeled 
$S_i $ with $i\in \Z/n\Z $ such that 
 $\sigma (S_i) = S_{i+1}$, where $\sigma $ is as in Theorem \ref{main2}. 

Then the simple $\Delta _s$-modules are 
$T_j$ with $j\in \Z/n' \Z $ can be labeled such that 
%the ${\cal B} \epsilon _s $-module
%$T_j$ is the direct sum of the $d$ simple modules 
$$(T_j) _{|{\cal B} \epsilon _s} = \bigoplus _{i\in \nu ^{-1}(cj)} S_i $$

The $\Delta _s$-lattices in the simple ${\cal A} \epsilon _s$-module 
form a chain 
$$\ldots \supset L_1\supset L_2 \supset \ldots \supset L_{n'} \supset pL_{1} =: L_{n'+1} \supset \ldots $$
where $L_j/L_{j+1} \cong T_j $ for $j=1,\ldots , n'$.
\end{rem}

\begin{rem}\label{headzp}
Theorem \ref{main2} also holds when the block ${\cal B}$ of $RG$
is replaced by the block $B$ of $\Z_pG$ from Theorem \ref{unique2}.
\end{rem}

\begin{rem}
Replacing $R$ by a ramified extension of $\Z_p$ in Remark \ref{headun}
and $a$ by the $\pi $-adic valuation of $p^a$ still yields a description 
of the head order of the non exceptional vertex $\Gamma _0$.
\end{rem}

\renewcommand{\arraystretch}{1}
\renewcommand{\baselinestretch}{1}
\large

\end{document}